\documentclass[letterpaper, 10 pt, conference]{ieeeconf}  
\IEEEoverridecommandlockouts  
\title{\LARGE \bf
Attack Detection in Dynamic Games with Quadratic Measurements
}

\usepackage{amsmath,amssymb,xcolor}
\usepackage{bm}
\usepackage{mathtools}
\usepackage{float}
\usepackage{stfloats} 
\usepackage{cite}

\DeclareMathOperator*{\argmin}{arg\,min}

\newtheorem{theorem}{Theorem}
\newtheorem{lemma}{Lemma}
\newtheorem{proposition}{Proposition}

\newtheorem{assumption}{Assumption}

\author{Muyan Jiang and Anil Aswani%
\thanks{*This material is based upon work supported by the National
Science Foundation under Grant DGE-2125913.}%
\thanks{The authors are with the department of Industrial Engineering and Operations Research, 
University of California, Berkeley, CA 94720, USA. 
{\tt\small \{muyan\_jiang, aaswani\}@berkeley.edu}}}

\begin{document}

\maketitle
\thispagestyle{empty}
\pagestyle{empty}

\begin{abstract}
This paper studies attack detection for discrete-time linear systems with stochastic process noise that produce both a vulnerable (i.e., attackable) linear measurement and a secured (i.e., unattackable) quadratic measurement. The motivating application of this model is a dynamic-game setting where the quadratic measurement is interpreted as a system-level utility or reward, and control inputs into the linear system are interpreted as control policies that, once applied, are known to all game participants and which steer the system towards a game-theoretic equilibrium (e.g., Nash equilibrium). To detect attacks on the linear channel, we develop a novel quadratic-utility-aware observer that leverages the secured quadratic output and enforces measurement consistency via a projection step. We establish three properties for this observer: feasibility of the true state, prox-regularity of the quadratic-constraint set, and a monotone error-reduction guarantee in the noise-free case. To detect adversarial manipulation, we compare linear and quadratic observer trajectories using a wild bootstrap maximum mean discrepancy (MMD) test that provides valid inference under temporal dependence. We validate our framework using numerical experiments of a pursuit–evasion game, where the quadratic observer preserves estimation accuracy under linear-sensor attacks, while the statistical test detects distributional divergence between the observers’ trajectories.
\end{abstract}

\section{Introduction}

Secure state estimation is critical for multi-agent systems in which multiple decision-makers coordinate actions from streamed sensor data~\cite{Mo2015,Pasqualetti2013}. There is extensive literature on resilience to false-data injection at the sensor/estimator level~\cite{Fawzi2014,Pasqualetti2013,Mo2015,hespanhol2017dynamic,pajic2016attack,Yang2020}, as well as anomaly detection methods that seek to identify unusual behaviors in data streams \cite{bhuyan2013network,blazquez2021review,pang2021deep}. However, there has been less work done on secure state estimation in multi-agent, dynamic games, which we distinguish from the literature that uses game-theoretic models of attacks on control system inputs and measurements \cite{LeBlanc2013,Zhu2015,Pawlick2019}.

This paper considers a discrete-time linear system with two types of measurements: In addition to the usual linear measurement, a single quadratic measurement is also made. We assume that the linear measurement can be attacked (i.e., corrupted by an adversary) while the quadratic measurement cannot be attacked. Though our model does not reference a multi-agent game, it is motivated by a game-theoretic setting where the quadratic measurement corresponds to a utility-function value or reward that is received by the entire system. The goal of this paper is two-fold: To develop an observer for quadratic measurements, and to develop a statistical testing framework to detect attacks on the linear measurements.

\subsection{Collusion Detection in Multi-Agent Games}
A closely related topic with increasing attention is detecting collusion in multi-agent games \cite{palshikar2008collusion,Mazrooei2013,hespanhol2020hypothesis,Bonjour2022,Greige2022}. One set of approaches that has been proposed to detect collusion is the use of statistical tests \cite{hespanhol2020hypothesis,Bonjour2022}. Another set of approaches leverage classical artificial intelligence (AI) \cite{palshikar2008collusion,Mazrooei2013,Greige2022}. For example, in large-scale team-based games, systems combining social networks with play metadata and unsupervised anomaly detection have been used to flag suspicious pairs \cite{Greige2022}. In repeated-game scenarios, model-agnostic regulatory audits that rely on statistical testing of observed behavior can reveal latent algorithmic collusion without requiring access to agents’ internal models or training procedures~\cite{hartline2024regulation}. Our work is related because it involves detecting undesired behavior in multi-agent games, but it differs in the type and model of undesired behavior.

\subsection{Observer Design for Quadratic Measurements}
Observer design for quadratic measurements is a less well-studied topic. One approach to observer design is to augment the state with derivatives or auxiliary variables of the quadratic output, which under certain conditions on the system convert the problem into an equivalent higher-dimensional linear one and enabling Kalman-like observers with convergence guarantees \cite{theodosis2021state}. Related efforts analyze control for linear–quadratic output systems, including stabilizability criteria \cite{Montenbruck2017}, and observability results for position estimation using only range or bearing data~\cite{Hamel2017}. These works address systems with only quadratic outputs and no adversarial interference. By contrast, we design and apply a novel observer for quadratic measurements to a system with a vulnerable linear channel, enforcing consistency with the quadratic measurement to yield an attack-resilient observer.

\subsection{Contributions and Outline}

We make two main contributions in this paper. The first is that we develop a novel state observer for quadratic measurements. The second is that we develop a statistical test that uses the quadratic measurement to  identify when the linear measurement is being attacked. 

Section \ref{sec:pf} presents the system model. Section \ref{sec:theory} defines our novel observer for quadratic measurements, and performs a theoretical analysis. Section \ref{sec:st} designs a statistical test for detecting attacks on the linear measurement. Section \ref{sec:exp} reports numerical experiments on a pursuit–evasion game, demonstrating detection of a sensor attack and maintenance of estimation accuracy under an attack.



\section{System Model}

This section presents the discrete-time linear system and its measurement model, and then it provides a game-theoretic interpretation of the quadratic measurements in this model. 

\label{sec:pf}

\subsection{Dynamics and Measurements}

Consider a linear system $x_{k+1} = Ax_k + Bu_k + w_k$ with stochastic process noise
$w_k \sim \mathcal{N}(0,\,Q)$, where $x_k\in\mathbb{R}^n$ is the state at time $k$, $u_k\in\mathbb{R}^m$ is the control input known to the estimator and treated as an exogenous input,
, $A\in\mathbb{R}^{n\times n}$ and $B\in\mathbb{R}^{n\times m}$ are known system matrices, and $w_k\in\mathbb{R}^n$ is zero‐mean Gaussian noise with covariance $Q\succ0$.

At each time step, two types of measurements are made: The linear measurements are $y_k = Cx_k + a_k + v_k$, with $v_k \sim \mathcal{N}(0,\,R)$,  where $C\in\mathbb{R}^{p\times n}$ is the measurement matrix, $v_k\in\mathbb{R}^p$ is Gaussian measurement noise, and $a_k\in\mathbb{R}^p$ is an unknown attack vector that may corrupt $y_k$. The quadratic measurements are $z_k = x_k^\top Vx_k$, with $V\in\mathbb{R}^{n\times n}$ and $V\succ0$, and this measurement cannot be manipulated by the attacker. 

\subsection{Game-Theoretic Interpretation}

Our motivation for studying the above model is the following game-theoretic interpretation: There are multiple agents subject to the dynamics, who select inputs according to an equilibrium policy (possibly nonlinear), and the realized input $u_k$ is known to all participants once applied. The linear measurements are susceptible to false-data injection \cite{Pasqualetti2013,Fawzi2014,pajic2016attack}, while the quadratic measurement is a game-theoretic, realized system utility or reward. The quadratic measurements may be a physical quantity (e.g., energy or Euclidean distance) obtained from local sensors and hence tamper-resistant \cite{vcapkun2002gps}.

\section{Observer for Quadratic Measurements}
\label{sec:theory}

Here, we present a novel observer for quadratic measurements. Then we theoretically analyze it.

\subsection{Observer Design}

Because the linear measurements are susceptible to attack, whereas the quadratic measurements are not, we use two observers: The first only uses linear measurements, and the second only uses quadratic measurements. Let $(\hat{x}^{L}_{k|k-1},\hat{x}^{L}_{k|k})$ and $(\hat{x}^{Q}_{k|k-1},\hat{x}^{Q}_{k|k})$ denote the prior and posterior state estimates of the linear and quadratic observers at time $k$, respectively.

For linear measurements, we use a Kalman filter~\cite{kalman1960}:
\[
\begin{aligned}
\textbf{Prediction:}&\quad
\hat{x}^{L}_{k|k-1}
= A\hat{x}^{L}_{k-1|k-1} + Bu_{k-1}, \\[4pt]
\textbf{Update:}&\quad
\hat{x}_{k|k}^{L}
= \hat{x}^{L}_{k|k-1}
+ L_k \bigl(y_k - C\hat{x}^{L}_{k|k-1}\bigr).
\end{aligned}
\]
where $P^{L}_{k|k-1}=A P^{L}_{k-1|k-1} A^\top+Q$ is predicted covariance, $P^{L}_{k|k}=(I-L_kC)P^{L}_{k|k-1}$ is updated covariance, and $L_k=P^{L}_{k|k-1}C^\top(C P^{L}_{k|k-1}C^\top+R)^{-1}$ is the Kalman gain. Since $a_k$ is unknown, the linear observer runs the nominal Kalman filter assuming $a_k\equiv 0$.
    
For quadratic measurements, we propose an extended‐Kalman–style observer, followed by a consistency projection:

\textbf{Prediction:} $\hat{x}_{k|k-1}^{Q} = A\hat{x}_{k-1|k-1}^{Q} + Bu_{k-1}$

\textbf{Extended Kalman Filter (EKF)‐Like Correction:}
    \begin{equation*}
      \tilde{x}_{k|k}
      = \hat{x}_{k|k-1}^{Q}
      + K_k \bigl(z_k - (\hat{x}_{k|k-1}^{Q})^\top V\hat{x}_{k|k-1}^{Q}\bigr),
    \end{equation*}
    where $H_k = \bigl(2V\hat{x}_{k|k-1}^{Q}\bigr)^\top$, $P^{Q}_{k|k}=(I-K_k H_k)P^{Q}_{k|k-1}$, $K_k = P^{Q}_{k|k-1}H_k^\top 
             (H_kP^{Q}_{k|k-1}H_k^\top + \eta I)^{-1}$, $\eta>0$ regularizes the gain, and $P^{Q}_{k|k-1}=A P^{Q}_{k-1|k-1} A^\top+Q$.\\
             
\textbf{Constrained Projection:}
\begin{equation*}      \hat{x}_{k|k}^{Q}
      = \argmin_{x \in \mathcal{F}_k}
        \|x - \tilde{x}_{k|k}\|_{{P^{Q}_{k|k}}^{-1}}^2,
\end{equation*}    
where $\mathcal{F}_k
      = \bigcap_{i=0}^{N}\{
        x : |H_{k-i}(A^{-i}x - \hat{x}_{k-i|k-i-1}^{Q}) - \tilde{z}_{k-i}|
        \le \delta_{k,i}(x)
      \},$ $\tilde{z}_{k-i}
      = z_{k-i} - (\hat{x}_{k-i|k-i-1}^{Q})^\top V\hat{x}_{k-i|k-i-1}^{Q}$, $\delta_{k,i}(x) 
      = \zeta + L\|A^{-i}x - \hat{x}_{k-i|k-i-1}^{Q}\|^2
$, $H_{k-i} = \bigl(2V\hat{x}_{k-i|k-i-1}^{Q}\bigr)^\top$ and 
      $L = \|V\|_2$. Here $\|\cdot\|_{P}^{2}$ denotes the squared weighted Euclidean norm for any symmetric $P\succeq 0$ and $N\ge 0$ denotes the number of past quadratic constraints enforced and $\zeta\ge 0$ is a tolerance term (taken as $\zeta=0$ in the noise-free analysis). The EKF‐like correction step treats the quadratic measurement \(z_k = x_k^\top V x_k\) as a nonlinear observation \(h(x)=x^\top V x\).  We linearize \(h\) around the prior \(\hat x_{k|k-1}^Q\) via its Jacobian \(H_k = (2V\hat{x}_{k|k-1}^Q)^\top\), and then apply a standard Kalman‐style update with gain \(K_k\).  

While this captures the local curvature of the quadratic sensor, it can drift when the linearization is poor. To counteract this, we project the corrected estimate \(\tilde x_{k|k}\) onto the feasible set \(\mathcal F_k\).  This set is defined by linearized measurement constraints from the current and past \(N\) steps, with adaptive bounds \(\delta_{k,i}(x)\) that account for the second‐order (linearization) error.  By solving the projection,
it returns the closest point, under the covariance‐weighted norm, to the unconstrained update, while remaining compatible with all secure quadratic measurements. This enhances robustness by anchoring the estimate to true system behavior, even in the presence of large innovations or attacked linear signals.


\subsection{Theoretical Error Bound}

Here, we analyze the noise-free case. Since the inputs $u_k$ are assumed to be known, without loss of generality we analyze our observer for the system: $x_{k+1} = Ax_k$ and $z_k = x_k^\top V x_k$, where $V$ is symmetric and positive definite.

We begin by noting that the absolute value constraint
\[
\bigl|H_{k-i}(A^{-i}x-\hat{x}^{Q}_{k-i|k-i-1})-\tilde{z}_{k-i}\bigr|\le \delta_{k,i}(x),
\]
 is equivalent to two inequalities. For $i=0,\dots,N$, define
\begin{align*}
\varphi_i^+(x) &= H_{k-i}(A^{-i}x-\hat{x}^{Q}_{k-i|k-i-1})-\tilde{z}_{k-i}-\delta_{k,i}(x),\\
\varphi_i^-(x) &= -H_{k-i}(A^{-i}x-\hat{x}^{Q}_{k-i|k-i-1})+\tilde{z}_{k-i}-\delta_{k,i}(x),
\end{align*}
so $\varphi_i^\pm(x)\le0$ encodes the same constraint. Since \(H_{k-i}\) and \(\tilde{z}_{k-i}\) are constants and \(\delta_{k,i}(x)\) is quadratic, each \(\varphi_i^\pm\) is \(C^2\).

Next we establish the prox-regularity of the feasible set $\mathcal{F}_k$ using the theory of amenable sets \cite{rockafellar1998variational}, by making some mild assumptions about constraint qualification.

\begin{assumption}[Invertibility]
\label{assumption:invertibility} $A\in\mathbb{R}^{n\times n}$ is invertible.
\end{assumption}
\begin{assumption}[Nondegeneracy]
\label{assumption:nondegeneracy}
For \(i=0,\dots,N\) and \(s\in\{+,-\}\), if \(\varphi_i^s(\bar{x})=0\) then
\(\nabla\varphi_i^s(\bar{x}) \neq 0\).
Equivalently, if \(y\in N_{(-\infty,0]}(\varphi_i^s(\bar{x}))\) and \(-\nabla\varphi_i^s(\bar{x})^*y=0\), then \(y=0\).
\end{assumption}

\begin{assumption}[Aggregated Constraint Qualification]
\label{assumption:constraint_qualification}
Define the stacked mapping
\begin{multline}
F(x) = \begin{bmatrix}
\varphi_0^+(x) & \varphi_0^-(x)  & \cdots &
\varphi_N^+(x) & \varphi_N^-(x)
\end{bmatrix}^\top .
\end{multline}
and let \(D=\prod_{j=1}^{2(N+1)}(-\infty,0]\), where $N_D(u)$ denotes the normal cone to $D$ at $u$. We assume that for \(\bar{x}\): if $y\in N_D\bigl(F(\bar{x})\bigr)$ and $\nabla F(\bar{x})^*y=0$, then $y=0$.
\end{assumption}

Unless stated otherwise, all results in this subsection hold under Assumptions~\ref{assumption:invertibility}--\ref{assumption:constraint_qualification}. We can formally define our feasible set as \(\mathcal{F}_k=\{x\in\mathbb{R}^n:F(x)\in D\}\). This formulation allows us to establish the main result:

\begin{proposition}[Prox-Regularity via Stacked Amenability]
Under Assumptions~\ref{assumption:invertibility}--\ref{assumption:constraint_qualification}, the set \(\mathcal{F}_k\) is strongly amenable at \(\bar{x}\) and, by \cite[Proposition 13.32]{rockafellar1998variational}, prox-regular at \(\bar{x}\).
\end{proposition}

\begin{proof}
Since each \(\varphi_i^\pm\) is \(C^2\) (due to its affine-plus-quadratic structure) and \(A\) is invertible by Assumption~\ref{assumption:invertibility}, the mapping \(F:\mathbb{R}^n\to\mathbb{R}^{2(N+1)}\) is \(C^2\). The set \(D=(-\infty,0]^{2(N+1)}\) is closed, convex, and polyhedral.

By \cite[Definition 10.23(b)]{rockafellar1998variational}, the representation \(\mathcal{F}_k=\{x\in \mathbb{R}^n:F(x)\in D\}\) establishes that \(\mathcal{F}_k\) is strongly amenable at \(\bar{x}\) provided the constraint qualification
\[
\text{if } y\in N_D(F(\bar{x})) \text{ and } \nabla F(\bar{x})^*y=0,\text{ then } y=0
\]
holds. Assumption~\ref{assumption:nondegeneracy} ensures that each active constraint \(\varphi_i^s\) is nondegenerate (i.e., \(\nabla\varphi_i^s(\bar{x})\neq 0\)), while Assumption~\ref{assumption:constraint_qualification} guarantees the aggregated constraint qualification for \(F\).

Therefore, by \cite[Proposition 13.32]{rockafellar1998variational}, the indicator function \(\delta_{\mathcal{F}_k}\) is prox-regular and subdifferentially continuous at \(\bar{x}\). Equivalently, the set \(\mathcal{F}_k\) is prox-regular at \(\bar{x}\).
\end{proof}

\begin{lemma}[Feasibility of State with Adaptive Bounds]
In the noise-free case, the true state $x_k$ belongs to the feasible set $\mathcal{F}_k$ when using the adaptive bounds $\delta_{k,i}(x) = \zeta + L\|A^{-i}x - \hat{x}_{k-i|k-i-1}^{Q}\|^2$ where $L = \|V\|_2$ and $\zeta = 0$ in the noise-free case.
\end{lemma}
\begin{proof}
For the true state $x_k$ to be in $\mathcal{F}_k$, it must satisfy:
$\bigl| \tilde{z}_{k-i} - H_{k-i}(A^{-i}x_k - \hat{x}_{k-i|k-i-1}^{Q}) \bigr|
\leq L\| A^{-i}x_k - \hat{x}_{k-i|k-i-1}^{Q} \|^2.$
From system dynamics (noise-free), $x_{k-i} = A^{-i}x_k$, so it suffices to verify
$\bigl|\tilde{z}_{k-i} - H_{k-i}(x_{k-i} - \hat{x}_{k-i|k-i-1}^{Q})\bigr|
\le L\|x_{k-i} - \hat{x}_{k-i|k-i-1}^{Q}\|^2.$
Let $e_{k-i|k-i-1} = x_{k-i} - \hat{x}_{k-i|k-i-1}^{Q}$. Substituting $H_{k-i}=(2V\hat{x}_{k-i|k-i-1}^{Q})^\top$ and $\tilde{z}_{k-i} = z_{k-i} - (\hat{x}_{k-i|k-i-1}^{Q})^\top V\hat{x}_{k-i|k-i-1}^{Q}$ gives
$\tilde{z}_{k-i} - 2(\hat{x}_{k-i|k-i-1}^{Q})^\top V e_{k-i|k-i-1} \leq L\|e_{k-i|k-i-1}\|^2.$ With $z_{k-i} = x_{k-i}^\top V x_{k-i}$ and $x_{k-i}=\hat{x}_{k-i|k-i-1}^{Q}+e_{k-i|k-i-1}$,
we have $\tilde{z}_{k-i}
= 2(\hat{x}_{k-i|k-i-1}^{Q})^\top V e_{k-i|k-i-1}
  + e_{k-i|k-i-1}^\top V e_{k-i|k-i-1}.$ Therefore,
$\bigl|\tilde{z}_{k-i} - H_{k-i}e_{k-i|k-i-1}\bigr|
= e_{k-i|k-i-1}^\top V e_{k-i|k-i-1}
\le \|V\|_2 \|e_{k-i|k-i-1}\|^2
= L\|e_{k-i|k-i-1}\|^2.$ Thus the $i$th constraint holds for $x_k$, and since $i$ is arbitrary, $x_k\in\mathcal{F}_k$.
\end{proof}

\begin{lemma}[Cross-Error Term Inequality]
\label{lem:ce-ineq}
Define the pre-projection error $\tilde{e}_{k+1|k+1} \triangleq \tilde{x}_{k+1|k+1} - x_{k+1}$ and the projection error $e^{obj}_{k+1} \triangleq \tilde{x}_{k+1|k+1} - \hat{x}_{k+1|k+1}$. Then, under prox-regularity of $\mathcal{F}_{k+1}$, for any $x_{k+1}\in \mathcal{F}_{k+1}$ (in particular for the true state) $\tilde{e}_{k+1|k+1}^{\top}(P^{Q}_{k+1|k+1})^{-1}e^{obj}_{k+1} \ge \|e^{obj}_{k+1}\|_{(P^{Q}_{k+1|k+1})^{-1}}^2$.
\end{lemma}

\begin{proof}
Since $\hat{x}_{k+1|k+1}$ is a local minimizer of 
\[
\hat{x}_{k+1|k+1} = \argmin_{x \in \mathcal{F}_{k+1}} \|x - \tilde{x}_{k+1|k+1}\|_{(P^{Q}_{k+1|k+1})^{-1}}^2,
\]
The first-order necessary optimality condition for constrained optimization requires $-\nabla f(\hat{x}_{k+1|k+1}) \in N_{\mathcal{F}_{k+1}}(\hat{x}_{k+1|k+1})$, where $N_{\mathcal{F}_{k+1}}(\hat{x}_{k+1|k+1})$ is the proximal normal cone to $\mathcal{F}_{k+1}$ at $\hat{x}_{k+1|k+1}$, and $\nabla f(x) = 2(P^{Q}_{k+1|k+1})^{-1}(x - \tilde{x}_{k+1|k+1})$. So
$2(P^{Q}_{k+1|k+1})^{-1}(\tilde{x}_{k+1|k+1} - \hat{x}_{k+1|k+1}) \in N_{\mathcal{F}_{k+1}}(\hat{x}_{k+1|k+1})$.

A key property of proximal normal cones for prox-regular sets is that for any $v \in N_{\mathcal{F}_{k+1}}(\hat{x}_{k+1|k+1})$ and any feasible point $x \in \mathcal{F}_{k+1}$, we have
$(x - \hat{x}_{k+1|k+1})^\top v \leq 0$ \cite{rockafellar1998variational}. Applying this to our case with $v = 2(P^{Q}_{k+1|k+1})^{-1}(\tilde{x}_{k+1|k+1} - \hat{x}_{k+1|k+1}) = 2(P^{Q}_{k+1|k+1})^{-1}e^{obj}_{k+1}$ and $x = x_{k+1}$, we get $(x_{k+1} - \hat{x}_{k+1|k+1})^T \cdot 2(P^{Q}_{k+1|k+1})^{-1}e^{obj}_{k+1} \leq 0$. Since $e^{obj}_{k+1} = \tilde{x}_{k+1|k+1} - \hat{x}_{k+1|k+1}$, substituting gives $x_{k+1} - \hat{x}_{k+1|k+1} = (x_{k+1} - \tilde{x}_{k+1|k+1}) + (\tilde{x}_{k+1|k+1} - \hat{x}_{k+1|k+1})
    = -\tilde{e}_{k+1|k+1} + e^{obj}_{k+1}$, which implies $(-\tilde{e}_{k+1|k+1} + e^{obj}_{k+1})^\top (P^{Q}_{k+1|k+1})^{-1} (-e^{obj}_{k+1}) \geq 0$. Expanding this gives $    \tilde{e}_{k+1|k+1}^\top (P^{Q}_{k+1|k+1})^{-1} e^{obj}_{k+1} - (e^{obj}_{k+1})^\top (P^{Q}_{k+1|k+1})^{-1} e^{obj}_{k+1} \geq 0$, which implies that 
$\tilde{e}_{k+1|k+1}^\top (P^{Q}_{k+1|k+1})^{-1} e^{obj}_{k+1} \geq \|e^{obj}_{k+1}\|_{(P^{Q}_{k+1|k+1})^{-1}}^2$.
\end{proof}

\begin{theorem}[Projection Error Bound]
Under the prox-regularity of $\mathcal{F}_{k+1}$, the projection step guarantees the post-projection error is bounded by the pre-projection error in the weighted norm $
\|e_{k+1|k+1}\|^2_{(P^{Q}_{k+1|k+1})^{-1}} \le \|\tilde{e}_{k+1|k+1}\|^2_{(P^{Q}_{k+1|k+1})^{-1}}$, where the post-projection error \footnote{Note that $e_{k+1|k+1}$ denotes the post-projection estimation error at time $k+1$, whereas the symbol $e_{k-i|k-i-1}$ in Lemma~1 is used locally to denote the prior estimation error at time $k-i$.
} is defined as $e_{k+1|k+1} \triangleq x_{k+1} - \hat{x}_{k+1|k+1}.$ 

\end{theorem}

\begin{proof}
We have $e_{k+1|k+1} = x_{k+1} - \hat{x}_{k+1|k+1} = (x_{k+1} - \tilde{x}_{k+1|k+1}) + (\tilde{x}_{k+1|k+1} - \hat{x}_{k+1|k+1}) = -\tilde{e}_{k+1|k+1} + e^{\mathrm{obj}}_{k+1}$. Thus, $\|e_{k+1|k+1}\|^2_{(P^{Q}_{k+1|k+1})^{-1}} = \|\tilde{e}_{k+1|k+1} - e^{\mathrm{obj}}_{k+1}\|^2_{(P^{Q}_{k+1|k+1})^{-1}} = \|\tilde{e}_{k+1|k+1}\|^2_{(P^{Q}_{k+1|k+1})^{-1}} + \|e^{\mathrm{obj}}_{k+1}\|^2_{(P^{Q}_{k+1|k+1})^{-1}} \quad - 2\,\tilde{e}_{k+1|k+1}^\top (P^{Q}_{k+1|k+1})^{-1} e^{\mathrm{obj}}_{k+1}$. By Lemma \ref{lem:ce-ineq}, we have $
\|e_{k+1|k+1}\|^2_{(P^{Q}_{k+1|k+1})^{-1}} 
\le \|\tilde{e}_{k+1|k+1}\|^2_{(P^{Q}_{k+1|k+1})^{-1}} - \|e^{\mathrm{obj}}_{k+1}\|^2_{(P^{Q}_{k+1|k+1})^{-1}} \le \|\tilde{e}_{k+1|k+1}\|^2_{(P^{Q}_{k+1|k+1})^{-1}}$.
\end{proof}

\section{Statistical Test for Attack Detection}

\label{sec:st}

This section develops a statistical test for detecting attacks on linear measurements using unattackable quadratic measurements. We test the null hypothesis $H_0$ that the state estimate distributions of the linear and quadratic observers coincide, against the alternative $H_1$ that an adversarial perturbation $a_k$ causes the linear observer’s estimates to deviate.

Since observer estimates are temporally dependent, standard permutation tests are invalid; we therefore adopt the wild bootstrap MMD test~\cite{chwialkowski2014wild}, which is designed for dependent data such as state trajectories. Let $X_k^{L} = \{\hat{x}_{1|1}^{L}, \ldots, \hat{x}_{k|k}^{L}\}$ and $X_k^{Q} = \{\hat{x}_{1|1}^{Q}, \ldots, \hat{x}_{k|k}^{Q}\}$. To quantify the discrepancy between these empirical distributions, we adopt the MMD equipped with an Radial Basis Function (RBF) kernel $
\phi(x,y) = \exp(-\|x - y\|^2/2\sigma^2)$. The empirical squared MMD is $\mathrm{MMD}^2(X_k^{L}, X_k^{Q})
= \tfrac{1}{k^2}\sum_{i,j=1}^k (\phi(\hat{x}_{i|i}^L,\hat{x}_{j|j}^L)
   + \phi(\hat{x}_{i|i}^Q,\hat{x}_{j|j}^Q) 
   - 2\phi(\hat{x}_{i|i}^L,\hat{x}_{j|j}^Q))$.

We combine the estimates into a single set $
Z_k = \{\hat{x}_{1|1}^{L},\dots,\hat{x}_{k|k}^{L}, \hat{x}_{1|1}^{Q},\dots,\hat{x}_{k|k}^{Q}\}$ with $2k$ total observations. Let $K\in\mathbb{R}^{2k\times 2k}$ with $K_{ij}=\phi(Z_i,Z_j)$ and define the centered kernel
$\tilde K=HKH$, where $H=I_{2k}-\frac{1}{2k}\mathbf{1}_{2k}\mathbf{1}_{2k}^\top$.
Let $\{v_i\}_{i=1}^{2k}$ be i.i.d.\ mean-zero, unit-variance random variables, and define
$\tilde K^v_{ij}=v_iv_j\tilde K_{ij}$.
The bootstrap statistic is
$\mathrm{MMD}_v=\frac{1}{2k}\sum_{i,j=1}^{2k}\tilde K^v_{ij}$.
This sum is a degenerate V-statistic and mimics the null distribution of $\mathrm{MMD}^2$ under temporal dependence. Repeating this procedure $B$ times yields $\{\mathrm{MMD}_v^{(b)}\}_{b=1}^B$, from which the critical value $\gamma_\alpha$ is obtained as the $(1-\alpha)$-quantile. Finally, the decision rule for attack detection becomes $\mathrm{MMD}^2(X_k^{L},X_k^{Q}) \underset{H_0}{\overset{H_1}{\gtrless}} \gamma_\alpha$.

\section{Numerical Experiments}
\label{sec:exp}
We conduct numerical experiments on a two-agent pursuit-evasion game governed by double integrator dynamics. Although our theoretical results assumed a noise-free regime, we include moderate Gaussian noise to demonstrate robustness beyond theoretical guarantees.

\subsection{Experimental Setup}
We consider a planar two-agent system with state vector $x_k \in \mathbb{R}^8$ at discrete time $k$, given by $x_k = [p_A, v_A, p_B, v_B]^\top$, where $p_A, p_B \in \mathbb{R}^2$ denote the positions and $v_A, v_B \in \mathbb{R}^2$ the velocities of the evader (Agent~A) and pursuer (Agent~B), respectively. The system evolves according to the discrete-time double integrator model $x_{k+1} = Ax_k + Bu_k + w_k$, where $A \in \mathbb{R}^{8 \times 8}$ and $B \in \mathbb{R}^{8 \times 4}$ are the state transition and input matrices, $u_k = [u_A, u_B]^\top$ is the control input, and $w_k$ is zero-mean Gaussian noise. The measurement model includes two channels: a vulnerable linear measurement $y_k = Cx_k + a_k + v_k$, where $C$ extracts the positions of both agents, $v_k \sim \mathcal{N}(0, R)$ is Gaussian noise, and $a_k$ is an adversarial attack vector; and a secure quadratic measurement $z_k = x_k^\top V x_k$, where $V \in \mathbb{R}^{8 \times 8}$ encodes the relative Euclidean distance.

The simulation parameters are as follows: time step $\Delta t = 0.1$\,s, simulation horizon $20$ steps, and process/measurements noise standard deviations all set to $0.005$. To assess robustness, we use randomized initial conditions drawn from Gaussian neighborhoods: the evader position $p_A(0)$ is sampled around $(0,0)$ with standard deviation $0.5$ in each axis, and the pursuer position $p_B(0)$ is sampled around $(2,2)$ with standard deviation $1.5$. Initial velocities have random directions (uniform over the unit circle) and magnitudes drawn from $\mathcal{N}(0.5,\,0.05^2)$ for the evader and $\mathcal{N}(0.2,\,0.05^2)$ for the pursuer, truncated below at $0.1$\,m/s. 

We run $M=100$ independent trials with the above randomized initializations. For each time step, we aggregate metrics across runs and report the mean and the standard error (SE). The attack is injected at discrete time $k=10$ with magnitude $\beta=7.0$ along the relative position direction.

\subsection{Control Policies}
\label{sec:control}


For discrete-time linear dynamics, optimal policies can be computed via  Hamilton--Jacobi--Bellman--Isaacs (HJBI) formulations \cite{Isaacs1965,Bardi1997,Falcone1994,tomlin2000game}, and capture conditions under full observability and sufficient control authority are well established \cite{Falcone1994}. For simplicity, we use heuristic control policies inspired by reachability-based strategies \cite{Chung2008} and observer-based estimation frameworks \cite{horak2017heuristic,oshman1999}. Control inputs are constrained component-wise by a saturation operator $[u]_{a_{\max}} \triangleq \max\{-a_{\max}, \min(u, a_{\max})\},\ a_{\max}=3\,\mathrm{m/s^{2}}$. This prevents physically unrealistic actuator demands.


\subsubsection*{Pursuer (Agent B, Leader)}

Agent~B uses perfect state knowledge to pursue an intercept point computed via one-step extrapolation with short-horizon interception timing:

\begin{enumerate}
    \item \textbf{Evader Prediction:} Predict the evader’s next position: $
    \tilde{p}_A(k+1) = p_A(k) + v_A(k)\Delta t$.
    \item \textbf{Intercept Calculation:} 
    Let $d_k \triangleq p_A(k)-p_B(k), \quad r_k \triangleq \|d_k\|_2$. If the evader is moving significantly (i.e., $\|v_A(k)\|_2 > 0.1\,\text{m/s}$), determine intercept time $t^\star$ by solving: $
    \|d_k + t(v_A(k)-v_B(k))\|_2^2 = (0.1\,r_k)^2,$,
    and set the intercept point as: $p_I = p_A(k)+v_A(k)\,t^\star$.
    Otherwise, default to the simple extrapolation: $p_I=\tilde{p}_A(k+1)$.
    \item \textbf{Desired Velocity:}
    The desired pursuit velocity combines range-dependent speed and near-range velocity matching: $
    v_B^{\mathrm{des}}(k)=
    s(r_k)\cdot(p_I - p_B(k))/\|p_I - p_B(k)\|_2+\beta(r_k)v_A(k)$,
    with: $s(r_k)= v_{\max,B}=2.5$ if $r_k>2$ and $s(r_k)= v_{\max,B}\cdot(0.5+0.25r_k)$ if $r_k\le 2$ and $\beta(r_k)=0.5$ if $r_k<1$ and $\beta(r_k)=0$ otherwise. The pursuer thus aggressively pursues at larger distances but smoothly transitions to cautious, velocity-matched intercept as the range closes, inspired by practical intercept strategies validated in~\cite{Chung2008}.

    \item \textbf{Control Law:}
    The control input for Agent~B is computed as $u_B(k)=[(v_B^{\mathrm{des}}(k)-v_B(k))/\Delta t]_{a_{\max}}$.
\end{enumerate}

\subsubsection*{Evader (Agent A, Follower)}

The evader relies exclusively on the observer estimate $\hat{x}_k$ and strategically evades by forecasting the pursuer’s short-term motion:

\begin{enumerate}
    \item \textbf{Pursuer Prediction:} Predict the pursuer’s next position from the estimate: $\tilde{p}_B(k+1)=\hat{p}_B(k)+\hat{v}_B(k)\Delta t$.
    \item \textbf{Escape Direction:} Compute the escape direction from the predicted pursuer position:$
    e_k=\hat{p}_A(k)-\tilde{p}_B(k+1)$.

    \item \textbf{Desired Velocity:}
    Maximize distance along the escape vector and add a minor velocity-matching perturbation to introduce unpredictability \cite{horak2017heuristic,Oshman2000} at longer distances: $    v_A^{\mathrm{des}}(k)
      = v_{\max,A}\cdot e_k/\|e_k\|_2
        + \gamma(\hat r_k)\,\hat{v}_B(k)$,
    with $v_{\max,A}=1.5$, 
    $\hat r_k = \|\hat p_A(k) - \hat p_B(k)\|_2$, and
    $\gamma(\hat r_k)=0.2$ if $\hat r_k>2$, otherwise $\gamma(\hat r_k)=0$.

    \item \textbf{Control Law:}
    The control input for Agent~A is similarly computed using estimated states: $u_A(k)=[(v_A^{\mathrm{des}}(k)-\hat{v}_A(k))/\Delta t]_{a_{\max}}$.
\end{enumerate}
The realized control inputs $u_k=[u_A(k);u_B(k)]$ are assumed to be known by both observers.

\subsection{Attack Scenario}
\label{sec:attack}
To evaluate detection and estimation robustness, we inject a \emph{relative position attack} on the linear channel at time $k=10$. The attack vector is constructed as $a_k = \beta \frac{p_B - p_A}{\|p_B - p_A\|}$, where $\beta = 7.0$ is the attack magnitude. We interpret $\beta$ as a distance bias magnitude in meters, injected along the relative position vector. This attack biases the perceived position of the pursuer, misleading the vulnerable observer.

\subsection{Experimental Results}

\begin{figure*}[t]
\centering
\includegraphics[width=0.9\textwidth]{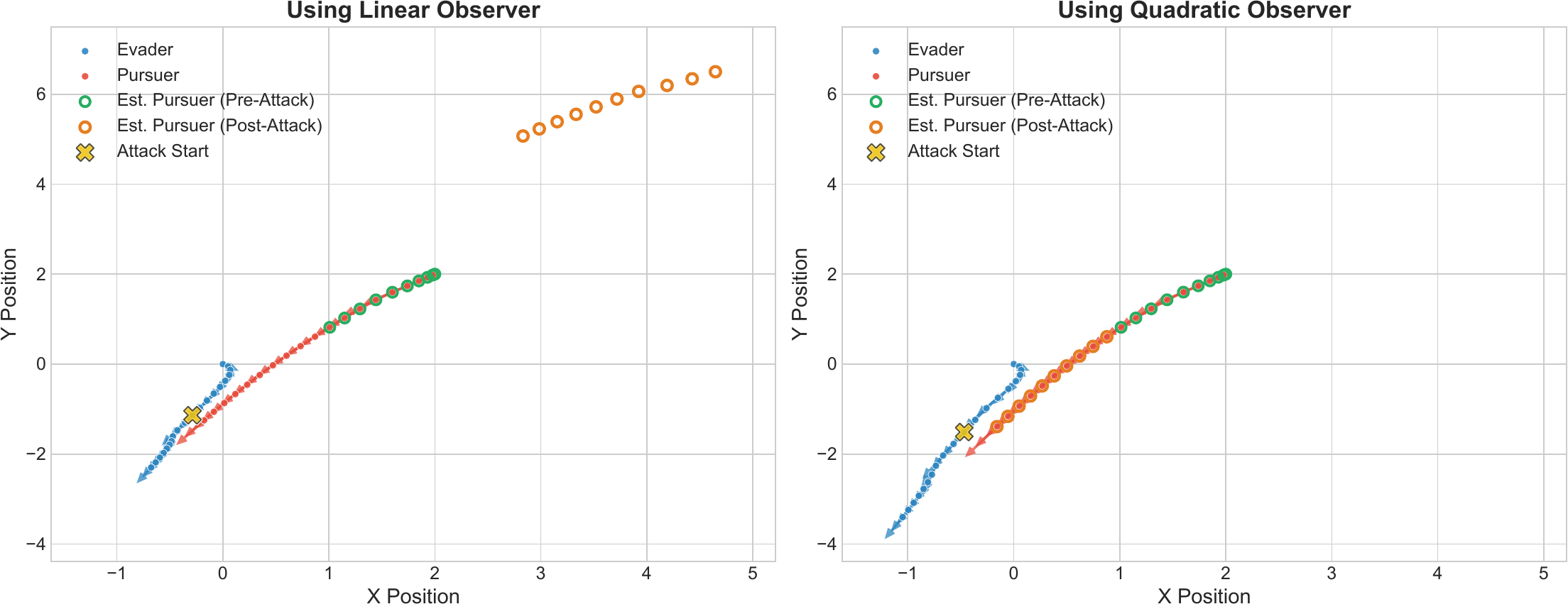}
\caption{Representative trial from the repeated experiments. Left: linear observer. Right: quadratic observer. True trajectories are shown for the evader (blue) and pursuer (red). The estimated pursuer trajectory is overlaid with hollow circles: green for pre–attack samples and orange for post–attack samples. Yellow ``X'' marks the attack onset. Faint lines trace the motion path and arrows indicate instantaneous velocity.}
\label{fig:trajectories}
\end{figure*}
\subsubsection{Trajectory Analysis}
\begin{figure}[h]
\centering
\includegraphics[width=0.95\linewidth]{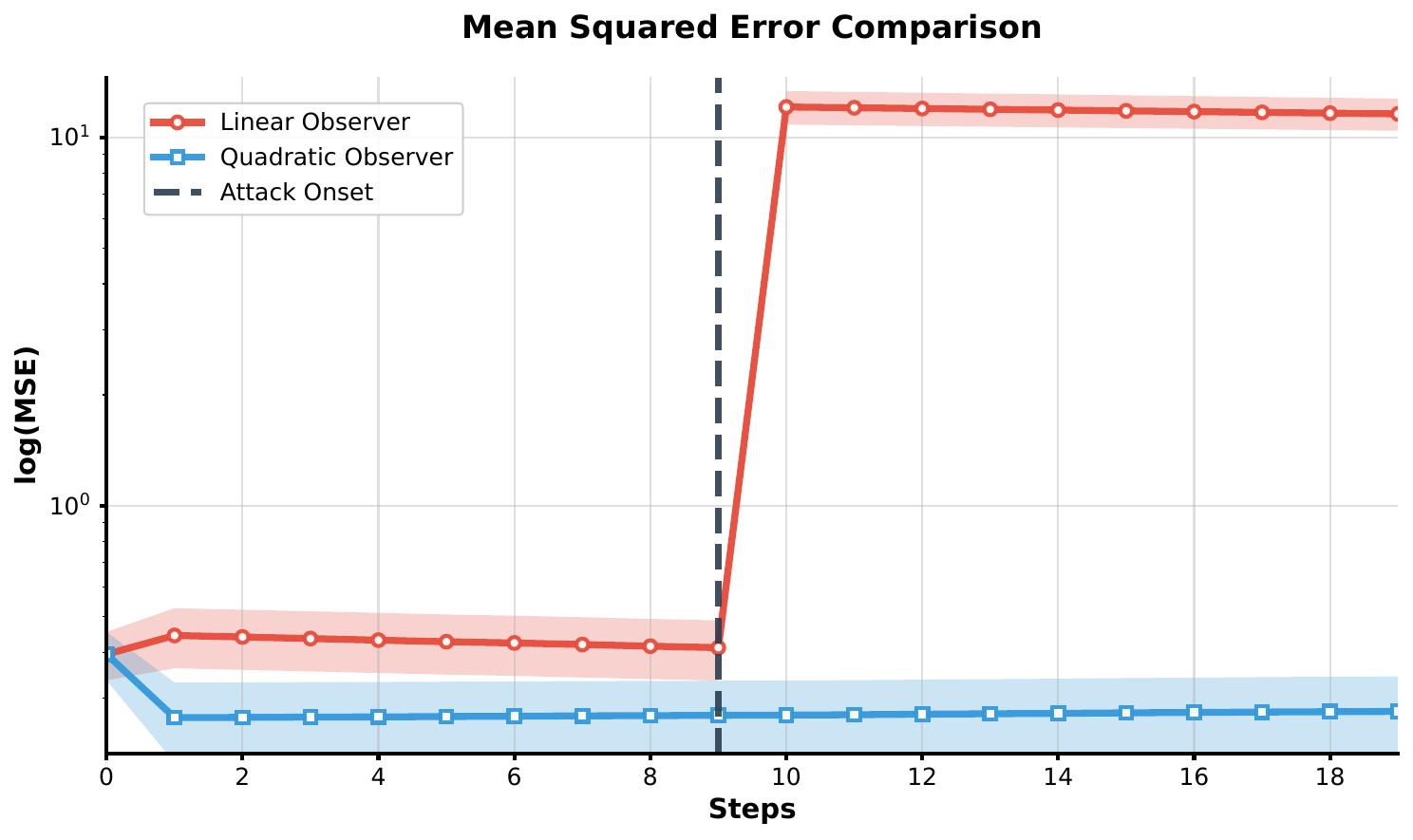}
\caption{Observer MSE over time aggregated across $M=100$ runs. Red: linear observer; blue: quadratic observer. Solid lines: mean MSE; shaded regions: $\pm$\,SE. Vertical dashed line indicates attack onset.}
\label{fig:mse}
\end{figure}

Fig.~\ref{fig:trajectories} shows a single representative trial drawn from the repeated-experiment protocol with randomized initial positions and velocities. Under the linear Kalman observer (left), the estimated pursuer trajectory (hollow orange circles, post-attack) departs from the red ground-truth path immediately after the attack marker (yellow ``X''). The drift appears as a systematic, directionally consistent bias that grows along the motion direction, yielding a spurious ``phantom'' pursuer that advances more slowly and farther from truth. In contrast, the quadratic observer (right) remains well aligned with the true pursuer trajectory both before and after the attack; the hollow green (pre-attack) and orange (post-attack) estimates closely overlay the red curve. Comparing the two panels over the same time horizon, the pursuer under the vulnerable linear observer appears to close the gap to the evader more than under the quadratic observer. Fig. \ref{fig:mse} summarizes the mean squared error (MSE) between true states and observer estimates across the repeated runs. Before the attack, both observers achieve comparable accuracy. After attack, the MSE of the linear observer increases markedly, whereas the quadratic observer maintains a low error by using the secure quadratic measurement.

\subsubsection{Attack Detection}
We use an RBF kernel (width via the median heuristic), and a wild bootstrap with Rademacher multipliers, $B=500$, $\alpha=0.05$. Online evaluation uses a sliding window $W$ equal to the pre-attack horizon; we declare detection at time $k$ if $\mathrm{MMD}^2_k>\hat\gamma_{\alpha,k}$. 

\begin{figure}[h]
\centering
\includegraphics[width=0.95\linewidth]{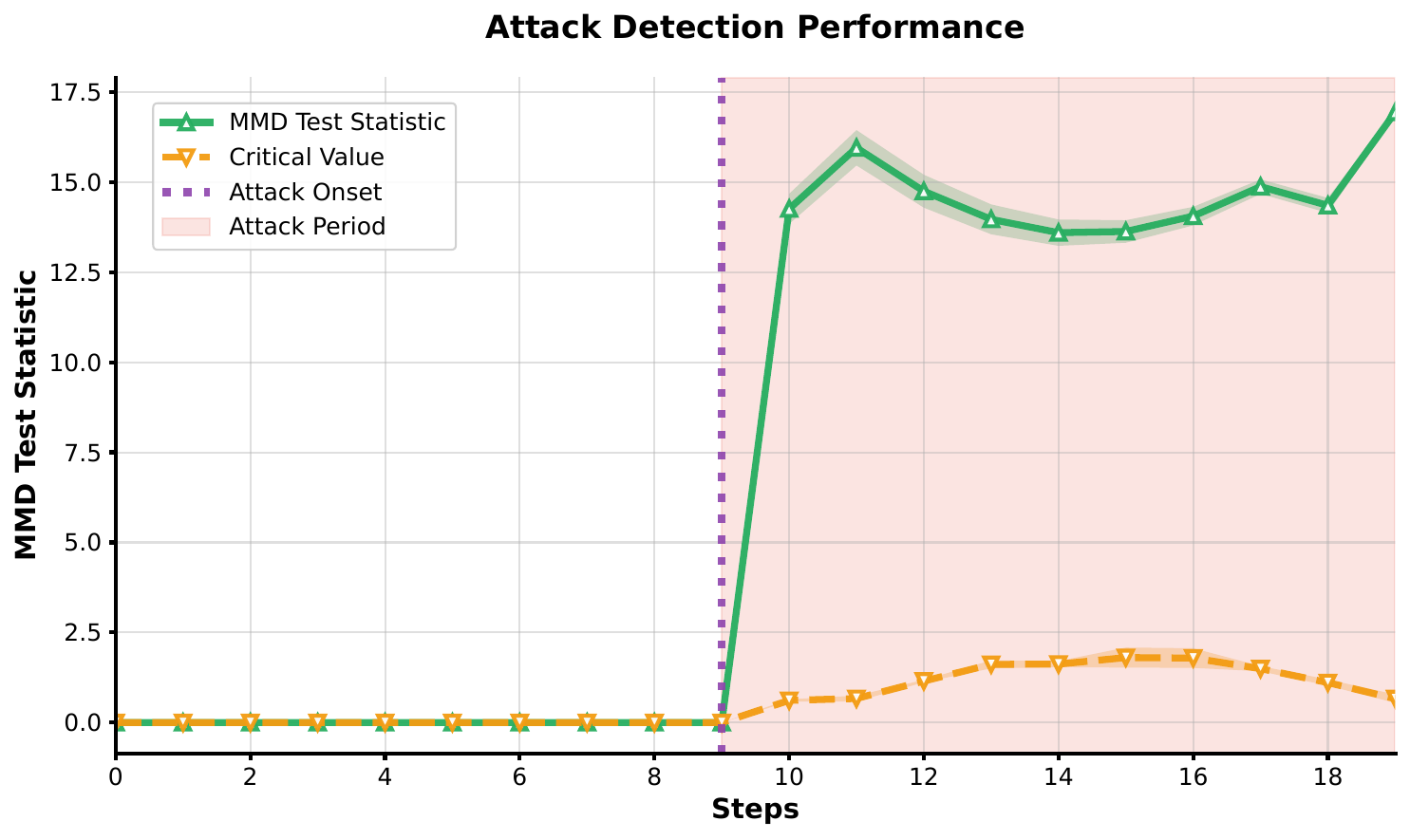}
\caption{Attack detection via wild bootstrap MMD over $M=100$ repeated experiments. Green: mean MMD test statistic; orange: mean critical value. Solid lines: mean; shaded bands: $\pm$\,SE. Vertical dotted line marks attack onset; red shaded region denotes the attack period.}
\label{fig:mmd}
\end{figure}

Fig.~\ref{fig:mmd} reports the \emph{aggregated} wild bootstrap MMD statistic across $M=100$ runs (mean $\pm$ SE) together with the corresponding mean critical value (dashed). Prior to the attack, the statistic remains below the threshold with no false positives on average. At the attack onset (vertical line), the mean MMD crosses the critical value with no delay, and the margin continues to widen thereafter, indicating a persistent distributional divergence between the drifted linear-observer trajectory and the stable quadratic-observer trajectory.

The results demonstrate that, with accurate initialization and low noise, the quadratic observer maintains robust state estimation in the presence of adversarial attacks, while the linear Kalman observer is significantly compromised. The MMD-based test provides prompt and reliable attack detection. These findings validate the theory and highlight the practical utility of the proposed approach for resilient state estimation and attack detection in dynamic games.

\section{Conclusion}

This work presented a robust framework for detecting adversarial sensor attacks in linear dynamical systems by combining a novel quadratic observer with a wild bootstrap MMD test. The quadratic observer leverages secure quadratic measurements to maintain reliable state estimates, while the wild bootstrap test detects distributional shifts under temporal dependence. Our theoretical analysis established error-monotonicity and prox-regularity properties of the proposed observer, and numerical experiments on a pursuit–evasion game demonstrated accurate estimation and prompt attack detection. Future work includes scaling the framework to larger multi-agent systems, incorporating adaptive thresholds for online testing, and extending the approach to nonlinear dynamics and broader classes of adversarial strategies.

\bibliographystyle{IEEEtran}
\bibliography{ref}

\end{document}